\documentclass[12pt]{article}

\usepackage{amsmath}
\usepackage{multirow} 
\usepackage{longtable}
\usepackage{latexsym,amsfonts}
\usepackage{graphicx}
\usepackage{subcaption}
\usepackage{pdflscape}
\usepackage{caption}
\usepackage{tikz-cd} 
\usepackage{tikz}
\usepackage{enumitem} 
\usepackage{authblk}
\newtheorem{thm}{Theorem}[section]

\newtheorem{prop}[thm]{Proposition}

\author[1]{Sanja Rukavina}
\author[2]{Vladimir D. Tonchev}
\affil[1]{Faculty of Mathematics, University of Rijeka, 51000 Rijeka, Croatia}
\affil[2]{ Department of Mathematical Sciences, Michigan Technological University, Houghton, MI 49931, USA
}
\title{Symmetric $2$-$(36,15,6)$ designs with an automorphism of order two}

\begin{document} 

\maketitle

\begin{abstract}

The parameters 2-$(36,15,6)$ are the smallest parameters of
symmetric designs for which a complete classification up to isomorphism is yet unknown. 
In \cite{BFW} Bouyukliev, Fack and Winne classified all 2-$(36,15,6)$ designs that admit
an automorphism of odd prime order, and gave a partial classification of such designs that 
admit an automorphism of order 2. In this paper, we give the classification of all symmetric 
2-$(36,15,6)$ designs that admit an automorphism of order two. It is shown that there are 
exactly $1 547 701$ nonisomorphic such designs,  $135 779$ of which are self-dual designs.
The ternary linear codes spanned by the incidence matrices of these designs are computed. 
Among these codes, there are near-extremal self-dual codes with previously unknown
weight distributions.
 
\end{abstract}

{\bf Mathematical subject classification (2020):} 05B05, 94B05 

\section {Introduction}
We assume familiarity with the basic facts and notions from the theory of combinatorial designs \cite{BJL,CRC, l, ton88}.

Symmetric designs with parameters $2$-$(36,15,6)$ can be obtained in many different ways, and it is not easy to determine the first occurrence of such a design. For example, in the first edition of Hall's  {\it Combinatorial Theory} \cite{Hall}, designs with parameters $2$-$(36,15,6)$ were listed as unkonwn, while in his 1964 work \cite{Mann}, Mann wrote about the connection between balanced incomplete block designs and abelian difference sets, and the existence of abelian difference sets with parameters $(36,15,6)$. These designs have attracted much attention, many connections have been made to other structures, and many such designs are presently known (see, e.g., \cite{BFW, BHSS, Sp}. A main open problem is the complete enumeration of such designs up to isomorphism. The problem of classifying all $2$-$(36,15,6)$ designs was approached by assuming an action of an automorphism group, leading to a partial classification. 
Bouyukliev, Fack and Winne \cite{BFW} gave a complete classification of
all  $2$-$(36,15,6)$ designs with an automorphism of odd prime order $p$,
and a partial classification for $p=2$. The partial classification from \cite{BFW} did not cover automorphisms of order two without fixed points, or involutions with more than $10$ fixed points.  
A classification of symmetric $2$-$(36,15,6)$ designs with an automorphism of order six was given in  \cite{ord9}, and all designs on which a cyclic automorphism group of order 4 acts standardly were constructed in \cite{crn}. In \cite{RT36}, $2$-$(36,15,6)$ designs with an automorphism of order two were used to investigate the existence of extremal self-dual ternary codes of length $36$. However, to the best of our knowledge, the enumeration of all symmetric $2$-$(36,15,6)$ designs with an automorphism of order two has been an open problem.
Many  examples of 2-$(36,15,6)$ designs having a trivial full automorphism group are known (see, for example, \cite{Sp}, \cite{Ton21}), but a complete classification of such designs is not known.\\
The links between Hadamard matrices of order 36, symmetric $2$-$(36,15,6)$ designs,
and extremal ternary self-dual codes of length 36 were studied in (\cite{RT36, Ton21}. 
In particular, it was proved in \cite{RT36} that, up to isomorphism, there is exactly one 
$2$-$(36,15,6)$ design 
with an an automorphism of order two whose incidence matrix spans an extremal ternary self-dual code of length $36$.  

The main goal of this paper is to give a detailed complete classification of symmetric $2$-$(36,15,6)$ designs with an automorphism of order two (Theorem \ref{main}). 
The ternary linear codes spanned by the incidence matrices of such designs are computed. 
Among these codes, there are near-extremal self-dual codes with previously unknown
weight distributions.

The paper is structured as follows. In the next section we consider all possible actions of an automorphism of order two on a symmetric $2$-$(36,15,6)$, and give in detail a classification of all such symmetric designs. In Section \ref{codes}, we use the obtained designs to   prove the existence of near-extremal self-dual ternary codes of length $36$ with $104$ or $120$ codewords of minimum weight.
 
\section{Classification of symmetric $2$-$(36,15,6)$  designs with an automorphism of order two}
\label{s2}

An automorphism of a symmetric design ${\mathcal D}$ is a permutation on the set of points of ${\mathcal D}$ that sends blocks to blocks. The action of an automorphism of ${\mathcal D}$ generates the same number of point and block orbits (see \cite[Theorem 3.3]{l}).
The set of all automorphisms of ${\mathcal D}$ forms its full automorphism group, denoted by $Aut({\mathcal D})$. A symmetric design is {\it self-dual} if it is isomorphic to its dual design\footnote{The {\it dual} design $D^T$ of a design $D$ with an incidence matrix $A$ is the design with the incidence matrix $A^T$.}.\\

A partial classification of symmetric $2$-$(36,15,6)$  designs with an automorphism of order two can be found in \cite{BFW}, and some information for the cases not treated there are given in \cite{RT36}. In both cases, the method for constructing orbit matrices with presumed action of an automorphism group, which are then indexed to construct designs was used (see, for example, \cite{ord9, jan}). We complete the classification using the same method. This work can be seen as a continuation of our research started in \cite{RT36}. Here we start with the orbit matrices considered in \cite{RT36}, which are available at \begin{verbatim}
 https://www.math.uniri.hr/~sanjar/structures/
 \end{verbatim}
In \cite{RT36} we were interested in a construction of extremal ternary self-dual codes of length $36$. Because of the large number of constructed designs, it was faster to check the codes corresponding to all designs than to eliminate the isomorphic copies of the designs first (see \cite[Note 2]{RT36}). Moreover, the only information given in \cite{RT36} about the case where an automorphism of order two acts on a symmetric $2$-$(36,15,6)$ design without fixed points is that such designs do not lead to self-dual ternary codes of length $36$. The same is true for symmetric $2$-$(36,15,6)$ designs acted upon by an automorphism of order two with six fixed points, for which it is falsely written that they do not exist. The number of constructed orbit matrices is therefore not given for all cases considered. Table \ref{tab3} contains complete information about the number of orbit matrices constructed in \cite{RT36}.

 \begin{table}[htpb!]
\begin{center} \begin{footnotesize}
\begin{tabular}{|c ||c| c | c|c|c|c|c| }
 \hline 

Number of fixed points& 0& 4 & 6 & 8 & 10 & 12 & 16 \\
 \hline
Number of orbit matrices&119909& 12991&2150& 670 & 56 & 311 & 83   \\
\hline
  
\end{tabular} \end{footnotesize}
 \caption{The number of orbit matrices for symmetric $2$-$(36,15,6)$ designs with an automorphism of order two}\label{tab3}
\end{center} 
\end{table}
  
In this paper we give details on the complete classification of all symmetric $2$-$(36,15,6)$ designs with an automorphism of order two, correcting some errors in the literature. In addition to our own computer programs, we use the System for Computational Discrete Algebra GAP \cite{GAP2022} and its package DESIGN \cite{Soi23} to determine the automorphism groups and representatives of the isomorphism classes of the constructed designs.\\

An automorphism of order two acting on $2$-$(36,15,6)$ design can act without fixed points or have $f$ fixed points, where $f \in \{4, 6, 8, 10, 12, 14, 16, 18 \}$. To our knowledge, a classification of symmetric $2$-$(36,15,6)$ designs for the cases $f \in \{0, 12, 14, 16, 18 \}$ cannot be found in the literature. In \cite{BFW} it is emphasized that for an action of an automorphism of order two, the most time-consuming part is the generation of the fixed parts (parts corresponding to orbits of length one) when  $f$ is large. Since there are no constraints for this part of the generation, the method used is too slow. We solved this problem in \cite{RT36} by reversing the order of the orbits and considering the rows corresponding to orbits of length two first. This imposes additional constraints on the construction of the fixed part in cases where all rows can be constructed for the block orbits of length two.\\

Once we have determined the possible orbit lengths distributions, we should consider the possible types for the rows of the orbit matrices corresponding to the orbits of certain lengths, i.e. the representatives of the possible types of block orbits  related to the given orbit lenghts distribution. For the cases where an automorphism of orbit two acts on the $2$-$(36,15,6)$ design with $14$ or $18$ fixed points, we will take a closer look at the corresponding row types to conclude that it is not possible to construct orbit matrices in these cases. For all other cases, our computational results are shown in the tables \ref{tab_inv_0f}-\ref{tab_inv_16f}, where we give the number of mutually non-isomorphic designs with the particular full automorphism group. We also give the number of self-dual designs for all cases. The final isomorphism test yields the classification shown in Table \ref{tab_inv_all}. The obtained results are compared with the previously known results.

\subsection {Symmetric $2$-$(36,15,6)$ designs admitting an automorphism of order two without fixed points} \label{sec-0f}

The possible action of an automorphism of order two on a symmetric $2$-$(36,15,6)$ design without fixed points was not mentioned in \cite{BFW}. Although this case was not specifically investigated, the results presented in \cite{ord9} show that such designs exist.\\

For this action of an automorphism of order two all orbits are of length two and all rows of the corresponding orbit matrices are of type
\begin{center}
2  2  2  1  1  1  1  1  1  1  1  1  0  0  0  0  0  0.
\end{center}
From $119909$ corresponding orbit matrices, we have constructed $13869$ mutually non-isomorphic $2$-$(36,15,6)$ designs (Table \ref{tab_inv_0f}). Among the constructed designs there are $1199$ self-dual designs.

\begin{table}[htpb!]
\begin{center} \begin{scriptsize}
\begin{tabular}{|c |c| c ||c|c|c|} 
 \hline   
 
The order & The structure &No. of&The order & The structure&No. of  \\
of $Aut\mathcal{(D)}$&of $Aut\mathcal{(D)}$&designs& of $Aut\mathcal{(D)}$&of $Aut\mathcal{(D)}$&designs\\
 \hline \hline

3888&$E_{81}:(S_4\times Z_2)$&1 &24 &$Z_3\times D_8$  &2 \\
\hline
1944&$(E_{27}:S_4)\times Z_3$&1&&$D_{24}$&2\\
\hline
648&$Z_3\times (E_{27}:D_8)$&1&&$Z_4\times S_3$&2\\
\hline
432&$((E_4:(E_9:Z_3)):Z_2):Z_2$&1&18&$Z_3\times S_3$&2\\
\hline
324&$E_9\times (E_9:Z_4)$&1&16&$(Z_4\times Z_2):Z_2$&15\\
\hline
&$Z_3\times (E_{27}:E_4)$&1&&$Z_2\times D_8$&3\\
\hline
240&$S_5\times Z_2$&1&12&$D_{12}$&11\\
\hline
216&$(E_4:(E_9:Z_3)):Z_2$&1&&$Z_6 \times Z_2$&1\\
\hline
144&$S_4\times S_3$&2&&$Z_{12}$&2\\
\hline
&$E_9:((Z_4\times Z_2):Z_2)$&1&&$A_4$&2\\
\hline
72&$Z_2\times S_3\times S_3$&2&8&$Z_4\times Z_2$&28\\
\hline
&$Z_3\times ((Z_6\times Z_2):Z_2)$&1&&$D_8$&22\\
\hline
&$(Z_3\times A_4):Z_2$&2&&$E_8$&18\\
\hline
48&$Z_2\times S_4$&2&6&$Z_6$&6\\
\hline
36&$Z_6\times S_3$&5&&$S_3$&1\\
\hline
&$S_3\times S_3$&1&4&$E_4$&693\\
\hline
&$Z_3\times (Z_3:Z_4)$&1&&$Z_4$&66\\
\hline
24&$S_4$&6&2&$Z_2$&12958\\
\hline
&$E_4\times S_3$&4&&&\\
\hline

\end{tabular} \end{scriptsize} 
 \caption{Symmetric $2$-$(36,15,6)$ designs with an automorphism of order two acting without fixed points} \label{tab_inv_0f}
\end{center} 
\end{table}

\subsection {Symmetric $2$-$(36,15,6)$ designs admitting an automorphism of order two acting with four fixed points} \label{sub-4f}

According to \cite{BFW} there are $170648$ symmetric $2$-$(36,15,6)$ designs on which an automorphism of order two acts with $4$ fixed points. After eliminating the isomorphic designs constructed from $12991$ orbit matrices for such an action, we also obtained $170648$ designs among which there are $17210$ self-dual designs. These include the design leading to an extremal ternary self-dual code of length $36$   equivalent to the Pless symmetry code $C(17)$ \cite{RT36}. Further details of the constructed designs are given in Table \ref{tab_inv_4f}.

\begin{table}[htpb!] 
\begin{center} \begin{scriptsize}
\begin{tabular}{|c |c| c ||c|c|c|}
 \hline   
 
The order & The structure &No. of&The order & The structure&No. of  \\
of $Aut\mathcal{(D)}$&of $Aut\mathcal{(D)}$&designs& of $Aut\mathcal{(D)}$&of $Aut\mathcal{(D)}$&designs\\
 \hline \hline

  51840&$O(5,3):Z_2$&1 &32&$E_4\times D_8$&12 \\
\hline
12096& $PSU(3,3):Z_2$&1 &&$(Z_8:Z_2):Z_2$&1 \\
\hline
3888&$E_{81}:(S_4\times Z_2)$&1 &24 & $Z_2\times A_4$&10 \\
\hline
1152&$Ex^+_{32}:(S_3\times S_3)$&1 &   & $S_4$&7 \\
\hline
648&$(S_3\times S_3\times S_3):Z_3$&1 &   & $E_4\times S_3$&6 \\
\hline
432&$S_3\times ((S_3\times S_3):Z_2)$&1 &&$SL(2,3)$&3\\
\hline
&$((E_4:(E_9:Z_3)):Z_2):Z_2$&1 & & $D_{24}$&3 \\
\hline
384&$Ex^+_{32}:D_{12}$&2 &&$(Z_6\times Z_2):Z_2$&2\\
\hline
360&$A_5\times S_3$&2 &&$Z_4\times S_3$&2\\
\hline
324&$(Z_3:(E_9:Z_3)):E_4$&1 & 18 & $Z_3\times S_3$&50 \\
\hline
240&$S_5\times Z_2$&1 & & $E_9 : Z_2$&18 \\
\hline
216&$(E_9:Q_8):Z_3$&1 & &$D_{18}$&1\\
\hline
162&$(E_9:Z_3):Z_2$&1 & 16&$D_{16}$&4\\
\hline
 144& $S_4\times S_3$&2 & & $Z_2\times D_8$&11  \\
 \hline
&$(Z_3\times SL(2,3)):Z_2$&2 & & $(Z_4\times Z_2):Z_2$&11 \\
\hline
&$E_9:((Z_4\times Z_2):Z_2)$&1 &&$Z_4\times E_4$&2\\
 \hline
108&$S_3\times (E_9:Z_2)$&1 &12&$A_4$&22\\
 \hline
96&$D_8\times A_4$&1 && $D_{12}$&65 \\
\hline
72& $Z_2\times S_3 \times S_3$&3 & & $Z_6 \times Z_2$&8  \\
\hline
&$Z_2\times (E_9:Z_4)$&2 &8 & $Z_4 \times Z_2$&118 \\
\hline
64&$D_8\times D_8$&2 & & $E_8$&96 \\
\hline
&$(E_8:E_4):Z_2$&1&   & $D_8$&134 \\
\hline
54&$E_9:Z_6$&6  &&$Z_8$&12\\
\hline
&$Z_9:Z_6$&2 &6 & $S_3$&315 \\
\hline
& $Z_3\times (E_9:Z_2)$&1 & & $Z_6$&458  \\
\hline
48& $Z_2\times S_4$&3 & 4&$Z_4$&1208 \\
\hline
&$(Z_3\times Q_8):Z_2$&1 && $E_4$&4405 \\
\hline
&$GL(2,3)$&1 & 2 & $Z_2$&163599 \\
\hline
  36 & $S_3\times S_3$&17& & & \\
  \hline
&$Z_2\times (E_9:Z_2)$&4 & & &\\
\hline

\end{tabular} \end{scriptsize} 
 \caption{Symmetric $2$-$(36,15,6)$ designs with an automorphism of order two acting with four fixed points} \label{tab_inv_4f}
\end{center} 
\end{table}

\subsection{Symmetric $2$-$(36,15,6)$ designs with an automorphism of order two fixing six points}

The existence of $135139$ symmetric $2$-$(36,15,6)$ designs with an automorphism of order two acting with six fixed points was proved in \cite{BFW}. Our construction yields the same number of mutually non-isomorphic designs shown in Table \ref{tab_inv_6f}. We obtain these designs from $2150$ orbit matrices constructed from three types of rows of the orbit matrices corresponding to block orbits of length one and four possible types of rows corresponding to block orbits of length two. The constructed designs contain $8325$ self-dual designs.

\begin{table}[htpb!] 
\begin{center} \begin{scriptsize}
\begin{tabular}{|c |c| c ||c|c|c|}
 \hline   
 
The order& The structure&No. of &The order& The structure&No. of \\
of $Aut\mathcal{(D)}$&of $Aut\mathcal{(D)}$&designs&of $Aut\mathcal{(D)}$&of $Aut\mathcal{(D)}$&designs\\
 \hline \hline

 3888&$E_{81}:(S_4\times Z_2)$&1 &  36 & $S_3\times S_3$&11  \\
\hline
1944& $(E_{27}:S_4)\times Z_3$&1 &  & $Z_6\times S_3$&7 \\
\hline
648& $Z_3\times (E_{27}:D_8)$&1 & 24 & $Z_3 \times D_8$&2  \\
\hline
 486& $(Z_3 \wr_{reg} Z_3) :S_3$&1 &   & $S_4$&9  \\
\hline
432 & $(((Z_6\times Z_6):Z_3):Z_2):Z_2$&1 &   & $E_4\times S_3$&4 \\
\hline
324& $(Z_3\times ((E_9:Z_3):Z_2)):Z_2$&1 && $D_{24}$&2 \\
\hline
& $Z_3\times ((Z_3\times (E_9:Z_2)):Z_2)$&1 && $Z_2\times A_4$&1 \\
\hline
240&$S_5\times Z_2$&1 & 18  & $E_9 : Z_2$&50 \\
\hline
216& $((Z_6\times Z_6):Z_3):Z_2$&1& & $Z_3\times S_3$&9  \\
\hline
162&$Z_3\times ((Z_9:Z_3):Z_2)$&3 & & $Z_6\times Z_3$&2 \\
\hline
&$Z_3\times ((E_9:Z_3):Z_2)$&1 &16 & $(Z_4\times Z_2):Z_2$&11 \\
\hline
&$E_9\times (E_9:Z_2)$&1 && $Z_2\times D_8$&3 \\
\hline
 144& $S_4\times S_3$&2 &12& $D_{12}$&41 \\
 \hline
 &$(Z_2\times (E_9:Z_4)):Z_2$&1 & & $Z_6 \times Z_2$&45 \\
 \hline
 108& $Z_3\times S_3 \times S_3$&2 &10&$D_{10}$&3 \\
 \hline
 &$S_3\times (E_9:Z_2)$&1 &8 & $Z_4 \times Z_2$&24 \\
 \hline
72& $Z_2\times S_3 \times S_3$&2 &     & $D_8$&53 \\
\hline
 & $Z_3\times ((Z_6 \times Z_2):Z_2)$&1 & & $E_8$&16 \\
 \hline
 & $(Z_3\times A_4):Z_2$&2 & 6 & $Z_6$&120 \\
\hline
54& $Z_3\times (E_9:Z_2)$&4 & & $S_3$&290 \\
\hline
 & $E_9\times S_3$&4 & 4 & $E_4$&2594 \\
\hline
48& $Z_2\times S_4$&2 & 2 & $Z_2$&131807 \\
\hline

\end{tabular} \end{scriptsize} 
 \caption{Symmetric $2$-$(36,15,6)$ designs with an automorphism of order two acting with six fixed points} \label{tab_inv_6f}
\end{center} 
\end{table}

\subsection{Symmetric $2$-$(36,15,6)$ designs with an automorphism of order two fixing eight points}

The number of symmetric $2$-$(36,15,6)$ designs constructed in \cite{BFW} with an automorphism of order two fixing eight points is $126817$. Starting from $670$ orbit matrices, we obtain $130414$ mutually non-isomorphic designs for this case. As in other cases, to eliminate isomorphic copies we have used
\begin{verbatim}
 BlockDesignIsomorphismClassRepresentatives();
\end{verbatim}
 in GAP. Among the constructed designs there are $16046$ self-dual designs. Further information is contained in Table \ref{tab_inv_8f}.\\
 
\begin{table}[htpb!] 
\begin{center} \begin{scriptsize}
\begin{tabular}{|c |c| c ||c|c|c|}
 \hline   
 
The order& The structure&No. of &The order& The structure&No. of \\
of $Aut\mathcal{(D)}$&of $Aut\mathcal{(D)}$&designs&of $Aut\mathcal{(D)}$&of $Aut\mathcal{(D)}$&designs\\
 \hline \hline

  51840&$O(5,3):Z_2$&1 &24&$Z_2\times A_4$ &5 \\
\hline
1152&$Ex^+_{32}:(S_3\times S_3)$ &1 &&$E_4\times S_3$&2\\
\hline
648&$(S_3\times S_3\times S_3):Z_3$&1& &$Z_3\times D_8$ &1\\
\hline
432&$S_3\times ((S_3\times S_3):Z_2)$&1&&$(Z_6\times Z_2):Z_2$&3\\
\hline
384&$Ex^+_{32}:D_{12}$ &2 &20&$Z_5:Z_4$&1\\
\hline
360&$A_5\times S_3$ &2 &18&$Z_3\times S_3$&18\\
\hline
320&$(E_{16}:Z_5):Z_4$&1&16&$E_{16}$&233\\
\hline
108&$E_{27}:Z_4$&2& &$Z_4\times E_4$ &2\\
\hline
96&$A_4\times D_8$&1&&$Z_4:Z_4$&2\\
\hline
&$(Z_2\times S_4):Z_2$&1& &$D_{16}$ &2\\
\hline
80&$E_{16}:Z_5$&2&&$(Z_4\times Z_2):Z_2$&14\\
\hline
72&$Z_2\times S_3\times S_3$&1&&$Z_2\times D_8$&7\\
\hline
  &$Z_2\times (E_9:Z_4)$&2&12&$A_4$ & 36\\
\hline
  64&$(E_8:Z_4):Z_2$&5&&$D_{12}$ & 28\\
\hline
  &$D_8\times D_8$&2&&$Z_6\times Z_2$ & 17\\
\hline
  54&$Z_3\times (E_9:Z_2)$&2&&$Z_3: Z_4$ & 24\\
\hline
  48&$E_4\times A_4$&5&10&$D_{10}$ & 3\\
\hline
  42&$Z_2\times (Z_7:Z_3)$&1&8&$D_8$ & 109\\
\hline
  36&$S_3\times S_3$&6&&$Z_4\times Z_2$ & 88\\
\hline
  &$Z_2\times (E_9:Z_2)$&2&&$E_8$ & 75\\
\hline
  &$Z_3\times A_4$&2&6&$Z_6$ & 373\\
\hline
32&$E_{16}:Z_2$&17&&$S_3$&14\\
\hline
  &$E_4\times D_8$&12&4&$E_4$ & 6107\\
\hline
  30&$Z_3\times D_{10}$&2&&$Z_4$ & 947\\
\hline
  &&&2&$Z_2$ & 122229\\
\hline
 
\end{tabular} \end{scriptsize} 
 \caption{Symmetric $2$-$(36,15,6)$ designs with an automorphism of order two acting with eight fixed points} \label{tab_inv_8f}
\end{center} 
\end{table}

\subsection{Symmetric $2$-$(36,15,6)$ designs with an automorphism of order two fixing ten points}

Our results for the number of symmetric $2$-$(36,15,6)$ designs that admit an action of an automorphism of order two with ten fixed points agree with those in \cite{BFW}. From $56$ orbit matrices constructed in \cite{RT36}, we obtained $64006$ mutually non-isomorphic designs, including $1282$ self-dual designs. Details about the automorphism groups of the constructed designs are given in Table \ref{tab_inv_10f}.

\begin{table}[htpb!] 
\begin{center} \begin{scriptsize}
\begin{tabular}{|c |c| c ||c|c|c|}
 \hline   
 
The order& The structure&No. of &The order& The structure&No. of \\
of $Aut\mathcal{(D)}$&of $Aut\mathcal{(D)}$&designs&of $Aut\mathcal{(D)}$&of $Aut\mathcal{(D)}$&designs\\
 \hline \hline

108&$Z_3\times S_3\times S_3 $&2&12&$D_{12} $ &2 \\
\hline
36&$Z_6\times S_3 $&2&&$Z_6\times Z_2 $ &44 \\
\hline
&$S_3\times S_3 $&2&6&$Z_6 $ &234 \\
\hline
18&$Z_3\times S_3 $&2&4&$E_4 $ &1500 \\
\hline
&&&2&$Z_2 $ &62218 \\
\hline

\end{tabular} \end{scriptsize} 
 \caption{Symmetric $2$-$(36,15,6)$ designs with an automorphism of order two acting with ten fixed points} \label{tab_inv_10f}
\end{center} 
\end{table}

\subsection{Symmetric $2$-$(36,15,6)$ designs with an automorphism of order two fixing $12$ points} \label{sub-12f}

Symmetric $2$-$(36,15,6)$ designs with an automorphism of order two fixing $12$ points have not been classified so far. Due to the large number of designs constructed from $311$ corresponding orbit matrices, the process of eliminating isomorphic copies was most challenging in this case. On the other hand, the results of \cite{RT36} show that only these designs and the designs for the case where an automorphism of order two fixes four points lead to self-dual ternary codes of length $36$, which makes this case of particular interest. After eliminating isomorphic copies, $981693$ mutually non-isomorphic designs are obtained, of which $79295$ are self-dual. The results of the analysis of the corresponding automorphism groups are shown in Table \ref{tab_inv_12f}.\\

\begin{table}[htpb!] 
\begin{center} \begin{scriptsize}
\begin{tabular}{|c |c| c ||c|c|c|}
 \hline   
 
The order& The structure&No. of &The order& The structure&No. of \\
of $Aut\mathcal{(D)}$&of $Aut\mathcal{(D)}$&designs&of $Aut\mathcal{(D)}$&of $Aut\mathcal{(D)}$&designs\\
 \hline \hline

  51840&$O(5,3):Z_2$&1 &24&$Z_2\times A_4$ &4\\
\hline
12096&$PSU(3,3):Z_2$&1&&$D_{24}$&3\\
\hline
1152&$Ex^+_{32}:(S_3\times S_3)$ &1 &&$(Z_6\times Z_2):Z_2$&1\\
\hline
432&$((S_3\times S_3):Z_2)\times S_3$&1&&$S_4$&2\\
\hline
384&$Ex^+_{32}:D_{12}$&2&&$E_4\times S_3$&2\\
\hline
320&$(E_{16}:Z_5):Z_4$&1&&$Z_3\times D_8$&1\\
\hline
144&$(Z_3\times SL(2,3)):Z_2$&2&18&$Z_3\times S_3$&7\\
\hline
108&$Z_3\times S_3\times S_3$&2&&$Z_6\times Z_3$&14\\
\hline
96&$D_8\times A_4$&1&16&$D_{16}$&4\\
\hline
&$(Z_2\times S_4):Z_2$&1&&$(Z_4\times Z_2):Z_2$&14\\
\hline
80&$E_{16}:Z_5$&2&&$E_{16}$&233\\
\hline
72&$Z_2\times S_3\times S_3$&1&&$Z_2\times D_8$&8\\
\hline
64&$(E_8:Z_4):Z_2$&5&&$Z_4:Z_4$&2\\
\hline
&$D_8\times D_8$&2&&$Z_4\times E_4$&2\\
\hline
&$(Z_8:E_4):Z_2$&1&12&$Z_{12}$&38\\
\hline
54&$E_9\times S_3$&2&&$D_{12}$&11\\
\hline
48&$E_4\times A_4$&5&&$Z_6\times Z_2$&57\\
\hline
&$Z_2\times S_4$&1&&$A_4$&2\\
\hline
&$(Z_3\times Q_8):Z_2$&1&8&$E_8$&80\\
\hline
&$GL(2,3)$&1&&$Z_4\times Z_2$&98\\
\hline
36&$S_3\times S_3$&3&&$D_8$&136\\
\hline
&$Z_2\times (E_9:Z_2)$&2&&$Z_8$&2\\
\hline
&$Z_6\times S_3$&2&6&$Z_6$&1187\\
\hline
32&$E_{16}:Z_2$&17&4&$E_4$&8657\\
\hline
&$E_4\times D_8$&12&&$Z_4$&797\\
\hline
&$(Z_8:Z_2):Z_2$&1&2&$Z_2$&970260\\
\hline

\end{tabular} \end{scriptsize} 
 \caption{Symmetric $2$-$(36,15,6)$ designs with an automorphism of order two acting with $12$ fixed points} \label{tab_inv_12f}
\end{center} 
\end{table}

\subsection{Symmetric $2$-$(36,15,6)$ designs with an automorphism of order two fixing $14$ points}

As mentioned above, the main obstacle in constructing orbit matrices for an automorphism of order two acting on symmetric $2$-$(36,15,6)$ designs with a large number of fixed points is the construction of the fixed part for such orbit matrices. Therefore, we first restrict our attention to the rows corresponding to block orbits of length two. In this case, an orbit matrix should have $11$ such rows, and there are only three possible types for these rows (the first $14$ digits in each type correspond to the point orbits of length one):
\begin{center}
1  1  1  1  1  1  0  0  0  0  0  0  0  0  1  1  1  1  1  1  1  1  1  0  0 \\
1  1  1  1  0  0  0  0  0  0  0  0  0  0  2  1  1  1  1  1  1  1  1  1  0 \\
1  1  0  0  0  0  0  0  0  0  0  0  0  0  2  2  1  1  1  1  1  1  1  1  1 \\
\end{center}

While trying to construct orbit matrices starting from fixed orbits takes too long and fails, this approach shows that it is possible to construct at most seven rows of an orbit matrix with row types corresponding to block orbits of length two. Since we must have $11$ rows corresponding to block orbits of length two, this means that an automorphism of order two cannot act on a symmetric $2$-$(36,15,6)$ design with $14$ fixed points.

\subsection{Symmetric $2$-$(36,15,6)$ designs with an automorphism of order two fixing $16$ points}

From $83$ orbit matrices for the action of an automorphism of order two on a symmetric $2$-$(36,15,6)$ design with $16$ fixed points constructed in \cite{RT36}, we obtain $67027$ mutually non-isomorphic designs. Among them, there are $20467$ self-dual designs. This is the first time that such designs are classified, and details about their automorphism groups are given in Table \ref{tab_inv_16f}.\\

\begin{table}[htpb!] 
\begin{center} \begin{scriptsize}
\begin{tabular}{|c |c| c ||c|c|c|}
 \hline   
 
The order& The structure &No. of&The order& The structure  &No. of  \\
of $Aut\mathcal{(D)}$&of $Aut\mathcal{(D)}$&designs&of $Aut\mathcal{(D)}$&of $Aut\mathcal{(D)}$&designs\\
 \hline \hline

 51840&$O(5,3):Z_2$&1 &24&$(Z_6\times Z_2):Z_2$&2 \\
\hline
3888&$E_{81}:(S_4\times Z_2)$&1 &&& \\
\hline
1152&$Ex^+_{32}:(S_3\times S_3)$&1& & $E_4\times S_3$&2   \\
\hline
648&$(S_3\times S_3\times S_3):Z_3$&1& 16& $Z_2\times D_8$&7  \\
\hline
432&$S_3\times ((S_3\times S_3):Z_2)$&1&&$Z_4\times E_4$&2  \\
\hline
384&$Ex^+_{32}:D_{12}$&2&& $(Z_4\times Z_2):Z_2$&4  \\
\hline
360&$A_5\times S_3$&2&12& $D_{12}$&22 \\
\hline
96&$D_8\times A_4$&1 & & $Z_6 \times Z_2$&18 \\
\hline
72& $Z_2\times S_3 \times S_3$&1& 10&$Z_{10}$&5  \\
\hline
&$Z_2\times (E_9:Z_4)$&2&8 & $Z_4 \times Z_2$&54 \\
\hline
  36 & $S_3\times S_3$&6& & $E_8$&67   \\
  \hline
64&$D_8\times D_8$&2&      & $D_8$&22 \\
\hline
36&$Z_2\times (E_9:Z_2)$&2 & 6 & $Z_6$&66 \\
\hline
32&$E_4\times D_8$&12 & 4 & $E_4$&1334  \\
\hline
24& $Z_2\times A_4$& 5& &$Z_4$&114\\
\hline
 & $Z_3 \times D_8$&1 & 2 & $Z_2$&65268\\
 \hline

\end{tabular} \end{scriptsize} 
 \caption{Symmetric $2$-$(36,15,6)$ designs with an automorphism of order two acting with $16$ fixed points} \label{tab_inv_16f}
\end{center} 
\end{table}

\subsection{Symmetric $2$-$(36,15,6)$ designs with an automorphism of order two fixing $18$ points}

There is only one type for the rows of the orbit matrix corresponding to the block orbits of length two for an action of an automorphism of order two  that fixes $18$ points on a symmetric $2$-$(36,15,6)$ design:

\begin{center}
1 1 1 1 1 1 0 0 0 0 0 0 0 0 0 0 0 0 1 1 1 1 1 1 1 1 1 \\
\end{center}

It is not difficult to see that it is not possible to combine two rows of this type to obtain the rows corresponding to the block orbits of length two of an orbit matrix for a symmetric $2$-$(36,15,6)$ design. Therefore, an automorphism of order two cannot act on a symmetric $2$-$(36,15,6)$ design with $18$ fixed points.

\subsection{Enumeration of symmetric $2$-$(36,15,6)$ designs with an automorphism of order two}

Here we summarise the results of the analysis of all possible cases for the action of an automorphism of order two on symmetric $2$-$(36,15,6)$ designs and give the main result of this paper. In the previous sections, we have eliminated isomorphic copies among the designs related to a particular action of an automorphism of order two. To obtain a classification of all symmetric $2$-$(36,15,6)$ designs with an automorphism of order two, we need to eliminate isomorphic designs that admit different actions of an automorphism of order two. If the full automorphism group of a design is of order two, then such a design admits only one action of an automorphism order two. However, if the full automorphism group is of composite order and contains different conjugacy classes of a subgroup of order two, it is possible that such a design admits different actions of an automorphism of order two (with respect to the number of fixed points). We have thus added the numbers of designs with the full automorphism group of order two for all possible actions considered in the previous sections and checked the remaining designs for isomorphism. Our analysis leads to a complete classification of symmetric $2$-$(36,15,6)$ designs admitting an automorphism of order two, which is summarised in the following theorem.

\begin{thm} \label{main}
An automorphism of order two acts on a symmetric $2$-$(36,15,6)$ design with $f$ fixed points, where $f \in \{0,4,6,8,10,12,16 \}$. Up to isomorphism, there are $1 547 701$ symmetric $2$-$(36,15,6)$ designs admitting an automorphism of order two with the full automorphism groups as given in Table  \ref{tab_inv_all}. Amog them there are $ 135779$ self-dual designs.
\end{thm}

\begin{table}[htpb!] 
\begin{center} \begin{scriptsize}
\begin{tabular}{|c |c| c ||c|c|c|}
 \hline   
 
The order& The structure&No. of &The order& The structure&No. of \\
of $Aut\mathcal{(D)}$&of $Aut\mathcal{(D)}$&designs&of $Aut\mathcal{(D)}$&of $Aut\mathcal{(D)}$&designs\\
 \hline \hline

  51840&$O(5,3):Z_2$&1 &48 &$E_4\times A_4$  &5\\
\hline
12096&$PSU(3,3):Z_2$&1&42&$Z_2\times (Z_7:Z_3)$&1\\
\hline
3888&$E_{81}:(S_4\times Z_2)$&1 &36&$S_3\times S_3$&19\\
\hline
1944& $(E_{27}:S_4)\times Z_3$&1 &&$Z_3\times A_4$&2\\
\hline
1152&$Ex^+_{32}:(S_3\times S_3)$&1&&$Z_6 \times S_3$&7\\
\hline
648&$(S_3\times S_3\times S_3):Z_3$&1&&$Z_2\times (E_9:Z_2)$&4\\
\hline
&$Z_3\times (E_{27}:D_8)$&1&&$Z_3\times (Z_3:Z_4)$&1\\
\hline
486&$(Z_3 \wr_{reg} Z_3) :S_3$&1&32&$E_4\times D_8$&12\\
\hline
432&$S_3\times ((S_3\times S_3):Z_2)$&1&&$E_{16}:Z_2$&17  \\
\hline
&$((E_4:(E_9:Z_3)):Z_2):Z_2$&1&&$(Z_8:Z_2):Z_2$&1\\
\hline
384&$Ex^+_{32}:D_{12}$&2&30&$Z_3\times D_{10}$&2\\
\hline
360&$S_3\times A_5$&2&24&$Z_2\times A_4$&10\\
\hline
324&$E_9\times (E_9:Z_4)$&1&&$S_4$&13\\
\hline
&$(Z_3\times (E_9:Z_3)):E_4$&1&&$SL(2,3)$&3\\
\hline
&$Z_3\times (E_{27}:E_4)$&1&&$D_{24}$&5\\
\hline
320&$(E_{16}:Z_5):Z_4$&1&&$Z_3\times D_8$&3\\
\hline
240&$Z_2\times S_5$&1&&$(Z_6\times Z_2):Z_2$&3\\
\hline
216&$(E_4:(E_9:Z_3)):Z_2$&1&&$E_4\times S_3$&6\\
\hline
&$(E_9:Q_8):Z_3$&1&&$Z_4\times S_3$&2\\
 \hline
162&$(E_{27}:Z_3):Z_2$&1&20&$Z_5:Z_4$&1\\
\hline
&$Z_3\times (Z_9:Z_6)$&3&18&$E_9:Z_2$&68\\
\hline
&$Z_3\times (E_9:Z_6)$&1&&$Z_6\times Z_3$&16\\
\hline
&$E_9\times (E_9:Z_2)$&1&&$D_{18}$&1\\
\hline
144&$S_3 \times S_4$&2&&$Z_3\times S_3$&88\\
\hline
&$(Z_3\times SL(2,3)):Z_2$&2&16&$D_{16}$&4\\
\hline
&$E_9:((Z_4\times Z_2):Z_2)$&1&&$(Z_4\times Z_2):Z_2$&25\\
\hline
108&$Z_3\times S_3\times S_3$&2&&$E_{16}$&233\\
\hline
&$E_{27}:Z_4$&2&&$Z_2\times D_8$&11\\
\hline
&$(E_9:Z_2)\times S_3$&1&&$Z_4:Z_4$&2\\
\hline
96&$(Z_2\times S_4):Z_2$&1&&$Z_4\times E_4$&2\\
\hline
&$D_8\times A_4$&1 & 12 & $A_4$&62 \\
\hline
80&$E_{16}:Z_5$&2&&$Z_{12}$&40\\
\hline
72&$Z_2\times S_3\times S_3$&3&&$Z_6\times Z_2$&70\\
\hline
&$Z_2\times (E_9:Z_4)$&2&&$D_{12}$&79\\
\hline
&$Z_3\times ((Z_6\times Z_2):Z_2)$&1&&$Z_3:Z_4$&24\\
\hline
&$(Z_3\times A_4):Z_2$&2&10&$D_{10}$&6\\
\hline
64&$D_8\times D_8$&2&&$Z_{10}$&5\\
\hline
&$(E_8:Z_4):Z_2$&5&8&$D_8$&235\\
\hline
&$(Z_8:E_4):Z_2$&1&&$Z_4\times Z_2$&180\\
\hline
54&$E_9:Z_6$&6&&$E_8$&108\\
\hline
&$Z_3\times (E_9:Z_2)$&7&&$Z_8$&14\\
\hline
&$E_9\times S_3$&6 & 6 & $Z_6$&2444 \\
\hline
&$Z_9:Z_6$& 2&  & $S_3$&620 \\
\hline
48&$Z_2\times S_4$ &3 & 4&$Z_4$&3132\\
\hline
&$(Z_3\times Q_8):Z_2$&1 & & $E_4$&11693  \\
\hline 
&$GL(2,3)$&1&2&$Z_2$&1528339\\
\hline

\end{tabular} \end{scriptsize} 
 \caption{Symmetric $2$-$(36,15,6)$ designs with an automorphism of order two} \label{tab_inv_all}
\end{center} 
\end{table}

For more detailed information on the structure of automorphism groups of orders greater than or equal to $486$, we refer the reader to \cite{ord9}.\\

When classifying combinatorial structures using algorithms and programs, which have been proven to provide accurate results in many previous cases, possible errors are mainly due to the human factor and errors that occur when entering, manipulating or analyzing a large amount of obtained data. In addition, when working with a very large data set, hardware-related errors are also possible, although these are very rare. Finally, the transcription of the final results is also subject to errors.

To verify the correctness of the results presented in this paper, we applied several techniques. First of all, we compared our results with known partial results. Our results are in complete agreement with the results published in \cite{BFW} for the cases where an automorphism of order two acts on a symmetric $2$-$(36,15,6)$ design with $4$, $6$ or $10$ fixed points, which gives us even more confidence in the procedures we have performed. Moreover, they agree with the results of \cite{ord9} for the vast majority of cases where the full automorphism group contains the subgroup that is isomorphic to the cyclic group of order six. For the few cases where our results do not match those in \cite{BFW,ord9}, we double-checked our results. These include the number of mutually non-isomorphic designs (we constructed more designs), the determinacy of the full automorphism group structure, and the pairing designs that are not self-dual. For the cases not previously discussed in the literature (an automorphism of order two acts on a symmetric $2$-$(36,15,6)$ design without fixed points or with more than $10$ points), we proceeded in the same way as for the confirmed previously known results. We performed the isomorphism test multiple times and with different groupings of designs, paired designs that are not self-dual, and checked the presented data multiple times to avoid errors in the presentation of the results.
 
\section{Near-extremal ternary self-dual codes of length $36$ related to the symmetric $2$-$(36,15,6)$ designs admitting an automorphism of order two} \label{codes}

We assume familiarity with the basic facts and concepts from the error-correcting codes \cite{BBFKKW, HP, Ton21}.

The sparsity of extremal self-dual ternary codes has recently aroused interest in near extremal ternary self-dual codes \cite{AR-HAR, HI, RT47}. In this section, we consider such codes that can be obtained from the symmetric $2$-$(36,15,6)$ designs with an automorphism of order two. In this part of our work we use the computer algebra system MAGMA \cite{magma}.\\

The minimum weight $d$ of a ternary self-dual code of length $n$ divisible by 12 satisfies the upper bound $d\le n/4 + 3$ \cite[9.3]{HP}. A self-dual ternary code of length $n$ divisible by 12 with a minimum weight $d$ is called {\it extremal} if $d=n/4+3$ \cite{HP}, and {\it near-extremal} if $d=n/4$ \cite{AR-HAR}.

Recently, it was proved in \cite{RT36} that there exists only one symmetric $2$-$(36,15,6)$ design with an automorphism of order two whose incidence matrix generates an extremal ternary self-dual code of length $36$. This design admits the action of an automorphism of order two with four fixed points and yields the only known extremal ternary self-dual code of length $36$, namely the Pless symmetry code $C(17)$ \cite{Pless69}.\\

The possible weight enumerators of near-extremal self-dual codes of length $36$ are determined in \cite{AR-HAR} and are given by
 $$1+\alpha y^9 +(42840-9\alpha)y^{12}+(1400256+36\alpha)y^{15}+(18452280-84\alpha)y^{18}$$
 \begin{equation}  \label{eq}
 +(90370368+126\alpha)y^{21}+(162663480-126\alpha)y^{24}+(97808480+84\alpha)y^{27} \end{equation}
$$+(16210656-36\alpha)y^{30}+(471240+9\alpha)y^{33}+(888-\alpha)y^{36}.$$

Moreover, the number of codewords of minimum weight in a ternary near-extremal self-dual code of length 36 is $A_9=8\beta$ for some $\beta$ in the range $1\le \beta \le 111$ \cite[page 1830]{AR-HAR}, and
 examples of near-extremal ternary self-dual codes of length $36$ are known for $A_9 \in \{8\beta | \beta \in \{9, 12, 14, 16, 17, . . . , 83, 85, 90, 93\} \}$ (see \cite{AR-HAR,HI}). The question of the existence of codes for other values of $A_9$ is still open.
Here we consider self-dual ternary codes arising from symmetric $2$-$(36,15,6)$ designs admitting an automorphism of order two in order to find examples of self-dual near-extremal ternary codes of length $36$ for values of $A_9$ for which the existence has not yet been proven.\\

As indicated in \cite[Table 3] {RT36}, an incidence matrix of a symmetric $2$-$(36,15,6)$ design with an automorphism of order two can generate a self-dual ternary code of length $36$ only if the corresponding automorphism of order two fixes four or twelve points. We consider all such designs from Section \ref{sub-4f} and Section \ref{sub-12f} and the results of our analysis are summarized in Table \ref{sd3}.

 \begin{table}[htpb!]
\begin{center} \begin{footnotesize}
\begin{tabular}{|c ||c| c | c|c|c|  }
 \hline 
&designs &self-dual&inequivalent& codes with& different \\ 
&&ternary codes& codes&$mw=9$&$A_9$\\
\hline
$4$ fixed points& 170648& 70774 & 2995  & 932&44   \\
 \hline
$12$ fixed points& $981693$& 391720&2363& 48&22    \\
\hline
  
\end{tabular} \end{footnotesize}
 \caption{Near-extremal self-dual ternary codes of length $36$ constructed from $2$-$(36,15,6)$ designs with an automorphism of order two}\label{sd3}
\end{center} 
\end{table}

The number of codewords with minimum weight in the obtained codes is given in Table \ref{A9}, where $i$ represents the number of fixed points for an action of an automorphism of order two on a symmetric $2$-$(36,15,6)$ design.

 \begin{table}[htpb!]
\begin{center} \begin{footnotesize}
\begin{tabular}{|c ||c|    }
 \hline 
i& $A_9$\\
\hline
4&$112, 128, 136, 152, 160, 168, 176, 184, 192, 200, 208, 216, 224, 232, 240, 248, 256, 264,$\\
&$272, 280, 288, 296, 304, 312, 320, 328, 336, 344, 352, 360, 368, 
376, 384, 392, 400, 408,$ \\
&$416, 424, 432, 440, 448, 464, 480, 544$    \\
 \hline
12&$104, 120, 136, 152, 168, 184, 200, 216, 232, 248, 264, 280, 296, 312, 328, 392, 424, 472,$ \\
&$488, 552, 584, 680$     \\
\hline
  
\end{tabular} \end{footnotesize}
 \caption{The number of codewords with minimum weight}\label{A9}
\end{center} 
\end{table}

Since near-extremal ternary self-dual codes with $A_9=8\cdot13=104$ and $A_9=8\cdot15=120$ were not constructed in \cite{AR-HAR, HI}, we have obtained first examples of such codes. Therefore, we present incidence matrices for one $2$-$(36,15,6)$ design that yield a self-dual near-extremal ternary code for each of these cases in the Appendix. Details on   analysis of obtained self-dual ternary $[36,18,9]$ codes with $A_9=104$ or $A_9=120$ are given in Table \ref{neq}.  \\ 

 \begin{table}[htpb!]
\begin{center} \begin{footnotesize}
\begin{tabular}{|c ||c| c | c|  }
 \hline 
&number of designs & inequivalent $[36,18,9]$& the order of $Aut\mathcal{(C)}$ \\
\hline
$A_9=104$& 12& 2 & 32, 48   \\
 \hline
$A_9=120$& 12& 7 & 24   \\
\hline
  
\end{tabular} \end{footnotesize}
 \caption{Self-dual ternary $[36,18,9]$ codes with $A_9=104$ or $A_9=120$}\label{neq}
\end{center} 
\end{table}

Our analysis, together with the results of \cite{AR-HAR, HI}, leads to the following conclusion.

\begin{prop}
There is a ternary near-extremal self-dual code of length 36 with weight enumerator (\ref{eq}) for $A_9=\alpha \in \{8\beta | \beta \in \{9, 12, . . . , 83, 85, 90, 93\} \}$.
\end{prop}

\bigskip
\noindent {\bf Acknowledgement} 

The first author is supported by {\rm C}roatian Science Foundation under the project 4571.

\section{Appendix}
The incidence matrix of a symmetric $2$-$(36,15,6)$ design that yield a near-extremal ternary self-dual code of length $36$ with $A_9=104$.
\begin{scriptsize}
\begin{verbatim}
  [ 1 1 1 1 1 1 1 0 0 0 0 0 1 1 1 1 1 1 1 1 0 0 0 0 0 0 0 0 0 0 0 0 0 0 0 0 ]
  [ 1 1 1 1 1 1 0 1 0 0 0 0 0 0 0 0 0 0 0 0 1 1 1 1 1 1 1 1 0 0 0 0 0 0 0 0 ]
  [ 1 1 1 0 0 0 1 1 1 1 0 0 1 1 0 0 0 0 0 0 1 1 0 0 0 0 0 0 1 1 1 1 0 0 0 0 ]
  [ 1 1 0 1 0 0 1 1 1 1 0 0 0 0 1 1 0 0 0 0 0 0 1 1 0 0 0 0 0 0 0 0 1 1 1 1 ]
  [ 1 1 0 0 0 0 0 0 0 0 1 0 1 1 1 1 0 0 0 0 0 0 0 0 1 1 1 1 1 1 0 0 0 0 1 1 ]
  [ 1 1 0 0 0 0 0 0 0 0 0 1 0 0 0 0 1 1 1 1 1 1 1 1 0 0 0 0 1 1 0 0 0 0 1 1 ]
  [ 1 0 1 1 0 0 1 1 0 0 1 1 0 0 0 0 1 1 0 0 0 0 0 0 1 1 0 0 1 1 0 0 1 1 0 0 ]
  [ 1 0 1 0 0 1 0 0 0 1 1 1 1 0 1 0 1 0 0 0 1 0 1 0 0 0 1 0 0 0 1 0 1 0 1 0 ]
  [ 1 0 1 0 0 1 0 0 0 1 1 1 0 1 0 1 0 1 0 0 0 1 0 1 0 0 0 1 0 0 0 1 0 1 0 1 ]
  [ 1 0 0 1 1 0 0 0 0 1 0 0 0 0 1 0 1 0 0 1 0 1 0 0 1 0 0 1 1 0 1 1 0 1 1 0 ]
  [ 1 0 0 1 1 0 0 0 0 1 0 0 0 0 0 1 0 1 1 0 1 0 0 0 0 1 1 0 0 1 1 1 1 0 0 1 ]
  [ 1 0 0 0 1 1 1 0 1 0 0 1 1 0 0 0 0 0 1 0 0 0 0 1 1 0 1 0 1 0 1 0 0 1 0 1 ]
  [ 1 0 0 0 1 1 1 0 1 0 0 1 0 1 0 0 0 0 0 1 0 0 1 0 0 1 0 1 0 1 0 1 1 0 1 0 ]
  [ 1 0 0 0 0 0 0 1 1 0 1 0 1 0 0 1 0 1 1 1 0 1 1 0 1 0 0 1 0 0 1 0 1 0 0 0 ]
  [ 1 0 0 0 0 0 0 1 1 0 1 0 0 1 1 0 1 0 1 1 1 0 0 1 0 1 1 0 0 0 0 1 0 1 0 0 ]
  [ 0 1 1 1 0 0 1 1 0 0 1 1 0 0 0 0 0 0 1 1 0 0 0 0 0 0 1 1 0 0 1 1 0 0 1 1 ]
  [ 0 1 1 0 0 1 0 0 1 0 0 0 0 0 1 0 0 1 0 1 1 0 0 0 0 1 0 1 1 0 1 0 1 1 0 1 ]
  [ 0 1 1 0 0 1 0 0 1 0 0 0 0 0 0 1 1 0 1 0 0 1 0 0 1 0 1 0 0 1 0 1 1 1 1 0 ]
  [ 0 1 0 1 1 0 0 0 1 0 1 1 1 0 0 1 1 0 0 0 1 0 0 1 0 0 0 1 1 0 0 1 1 0 0 0 ]
  [ 0 1 0 1 1 0 0 0 1 0 1 1 0 1 1 0 0 1 0 0 0 1 1 0 0 0 1 0 0 1 1 0 0 1 0 0 ]
  [ 0 1 0 0 1 1 0 1 0 1 1 0 1 0 0 0 0 1 1 0 0 0 1 0 0 1 0 0 1 0 0 1 0 1 1 0 ]
  [ 0 1 0 0 1 1 0 1 0 1 1 0 0 1 0 0 1 0 0 1 0 0 0 1 1 0 0 0 0 1 1 0 1 0 0 1 ]
  [ 0 1 0 0 0 0 1 0 0 1 0 1 1 0 1 0 0 1 0 1 0 1 0 1 1 1 1 0 0 0 0 1 1 0 0 0 ]
  [ 0 1 0 0 0 0 1 0 0 1 0 1 0 1 0 1 1 0 1 0 1 0 1 0 1 1 0 1 0 0 1 0 0 1 0 0 ]
  [ 0 0 1 1 0 0 0 0 1 0 0 0 1 1 0 0 1 1 0 0 0 0 1 1 1 1 0 0 0 0 1 1 0 0 1 1 ]
  [ 0 0 1 1 0 0 0 0 0 1 0 0 1 1 0 0 0 0 1 1 0 0 1 1 0 0 1 1 1 1 0 0 1 1 0 0 ]
  [ 0 0 1 0 1 0 1 0 1 1 1 0 0 0 1 0 0 1 1 0 1 0 0 1 1 0 0 1 0 1 0 0 0 0 1 0 ]
  [ 0 0 1 0 1 0 1 0 1 1 1 0 0 0 0 1 1 0 0 1 0 1 1 0 0 1 1 0 1 0 0 0 0 0 0 1 ]
  [ 0 0 1 0 1 0 0 1 0 0 0 1 1 0 1 1 0 0 0 1 1 0 1 0 1 0 0 0 0 1 0 1 0 1 0 1 ]
  [ 0 0 1 0 1 0 0 1 0 0 0 1 0 1 1 1 0 0 1 0 0 1 0 1 0 1 0 0 1 0 1 0 1 0 1 0 ]
  [ 0 0 0 1 0 1 1 0 0 0 1 0 1 0 0 1 0 0 0 1 1 1 0 1 0 1 0 0 0 1 1 0 0 1 1 0 ]
  [ 0 0 0 1 0 1 1 0 0 0 1 0 0 1 1 0 0 0 1 0 1 1 1 0 1 0 0 0 1 0 0 1 1 0 0 1 ]
  [ 0 0 0 1 0 1 0 1 1 1 0 1 1 0 1 0 1 0 1 0 0 1 0 0 0 1 0 1 0 1 0 0 0 0 0 1 ]
  [ 0 0 0 1 0 1 0 1 1 1 0 1 0 1 0 1 0 1 0 1 1 0 0 0 1 0 1 0 1 0 0 0 0 0 1 0 ]
  [ 0 0 0 0 1 0 1 1 0 0 0 0 1 1 0 0 1 1 0 0 1 1 0 0 0 0 1 1 0 0 0 0 1 1 1 1 ]
  [ 0 0 0 0 0 1 1 1 0 0 0 0 0 0 1 1 1 1 0 0 0 0 1 1 0 0 1 1 1 1 1 1 0 0 0 0 ]
\end{verbatim}
\end{scriptsize}

The incidence matrix of a symmetric $2$-$(36,15,6)$ design that yield a near-extremal ternary self-dual code of length $36$ with $A_9=120$.
\begin{scriptsize}
\begin{verbatim}
    [ 1 1 1 1 1 1 1 0 0 0 0 0 1 1 1 1 1 1 1 1 0 0 0 0 0 0 0 0 0 0 0 0 0 0 0 0 ]
    [ 1 1 1 1 0 0 0 1 1 1 0 0 1 1 0 0 0 0 0 0 1 1 1 1 1 1 0 0 0 0 0 0 0 0 0 0 ]
    [ 1 1 1 0 1 0 0 1 0 0 1 1 1 1 0 0 0 0 0 0 0 0 0 0 0 0 1 1 1 1 1 1 0 0 0 0 ]
    [ 1 1 1 0 0 0 1 0 0 0 0 0 0 0 1 0 0 0 1 0 1 0 1 0 1 0 1 0 1 0 1 0 1 0 1 1 ]
    [ 1 1 1 0 0 0 1 0 0 0 0 0 0 0 0 1 0 0 0 1 0 1 0 1 0 1 0 1 0 1 0 1 0 1 1 1 ]
    [ 1 1 0 1 0 1 0 0 1 0 1 1 0 0 1 1 0 0 0 0 1 1 0 0 0 0 1 1 0 0 0 0 1 1 0 0 ]
    [ 1 0 1 0 1 1 0 0 1 1 1 0 0 0 0 0 1 1 0 0 1 1 0 0 0 0 0 0 1 1 0 0 0 0 1 1 ]
    [ 1 0 0 1 1 1 0 1 0 1 0 1 0 0 1 1 0 0 0 0 0 0 1 1 0 0 0 0 0 0 1 1 0 0 1 1 ]
    [ 1 0 0 1 0 1 0 0 0 0 1 0 1 0 0 0 1 0 1 0 0 0 1 0 1 1 0 1 0 1 1 0 0 1 1 0 ]
    [ 1 0 0 1 0 1 0 0 0 0 1 0 0 1 0 0 0 1 0 1 0 0 0 1 1 1 1 0 1 0 0 1 1 0 0 1 ]
    [ 1 0 0 0 1 0 0 1 1 0 0 0 0 0 1 0 0 1 1 1 0 1 1 0 1 0 0 1 1 0 0 1 0 1 0 0 ]
    [ 1 0 0 0 1 0 0 1 1 0 0 0 0 0 0 1 1 0 1 1 1 0 0 1 0 1 1 0 0 1 1 0 1 0 0 0 ]
    [ 1 0 0 0 0 0 1 1 1 1 1 1 1 1 0 0 0 0 1 1 0 0 0 0 0 0 0 0 0 0 0 0 1 1 1 1 ]
    [ 1 0 0 0 0 0 1 0 0 1 0 1 1 0 0 1 1 1 0 0 1 0 1 0 0 1 0 1 1 0 0 1 1 0 0 0 ]
    [ 1 0 0 0 0 0 1 0 0 1 0 1 0 1 1 0 1 1 0 0 0 1 0 1 1 0 1 0 0 1 1 0 0 1 0 0 ]
    [ 0 1 1 0 0 1 1 1 1 0 1 0 0 0 0 0 1 1 0 0 0 0 1 1 0 0 0 0 0 0 1 1 1 1 0 0 ]
    [ 0 1 0 1 1 0 1 1 1 0 0 1 0 0 0 0 1 1 0 0 0 0 0 0 1 1 1 1 0 0 0 0 0 0 1 1 ]
    [ 0 1 0 1 1 0 0 0 0 1 0 0 1 0 0 0 1 0 1 0 0 1 0 1 0 0 1 0 1 0 0 1 1 1 1 0 ]
    [ 0 1 0 1 1 0 0 0 0 1 0 0 0 1 0 0 0 1 0 1 1 0 1 0 0 0 0 1 0 1 1 0 1 1 0 1 ]
    [ 0 1 0 0 1 1 1 0 0 1 1 1 0 0 0 0 0 0 1 1 1 1 0 0 1 1 0 0 0 0 1 1 0 0 0 0 ]
    [ 0 1 0 0 0 1 0 1 0 0 0 1 1 0 0 1 0 1 1 0 1 0 0 1 1 0 0 0 1 1 0 0 0 1 0 1 ]
    [ 0 1 0 0 0 1 0 1 0 0 0 1 0 1 1 0 1 0 0 1 0 1 1 0 0 1 0 0 1 1 0 0 1 0 1 0 ]
    [ 0 1 0 0 0 0 0 0 1 1 1 0 1 0 1 1 1 0 0 1 0 0 1 0 1 0 1 0 0 1 0 1 0 0 0 1 ]
    [ 0 1 0 0 0 0 0 0 1 1 1 0 0 1 1 1 0 1 1 0 0 0 0 1 0 1 0 1 1 0 1 0 0 0 1 0 ]
    [ 0 0 1 1 1 0 1 1 0 1 1 0 0 0 1 1 0 0 0 0 0 0 0 0 1 1 0 0 1 1 0 0 1 1 0 0 ]
    [ 0 0 1 1 0 1 1 0 1 1 0 1 0 0 0 0 0 0 1 1 0 0 1 1 0 0 1 1 1 1 0 0 0 0 0 0 ]
    [ 0 0 1 1 0 0 0 0 1 0 0 1 1 0 1 0 0 1 1 0 0 1 0 0 0 1 0 0 0 1 1 1 1 0 0 1 ]
    [ 0 0 1 1 0 0 0 0 1 0 0 1 0 1 0 1 1 0 0 1 1 0 0 0 1 0 0 0 1 0 1 1 0 1 1 0 ]
    [ 0 0 1 0 1 0 0 0 0 0 1 1 1 0 1 0 0 1 0 1 1 0 1 1 0 1 1 0 0 0 0 0 0 1 1 0 ]
    [ 0 0 1 0 1 0 0 0 0 0 1 1 0 1 0 1 1 0 1 0 0 1 1 1 1 0 0 1 0 0 0 0 1 0 0 1 ]
    [ 0 0 1 0 0 1 0 1 0 1 0 0 1 0 0 1 0 1 0 1 0 1 0 0 1 0 1 1 0 0 1 0 1 0 1 0 ]
    [ 0 0 1 0 0 1 0 1 0 1 0 0 0 1 1 0 1 0 1 0 1 0 0 0 0 1 1 1 0 0 0 1 0 1 0 1 ]
    [ 0 0 0 1 0 0 1 1 0 0 1 0 1 0 1 0 1 0 0 1 1 1 0 1 0 0 0 1 1 0 1 0 0 0 0 1 ]
    [ 0 0 0 1 0 0 1 1 0 0 1 0 0 1 0 1 0 1 1 0 1 1 1 0 0 0 1 0 0 1 0 1 0 0 1 0 ]
    [ 0 0 0 0 1 1 1 0 1 0 0 0 1 1 1 0 0 0 0 0 1 0 0 1 1 0 0 1 0 1 0 1 1 0 1 0 ]
    [ 0 0 0 0 1 1 1 0 1 0 0 0 1 1 0 1 0 0 0 0 0 1 1 0 0 1 1 0 1 0 1 0 0 1 0 1 ]
\end{verbatim}
\end{scriptsize}

\end{document}